\documentclass[12pt]{article}
\usepackage{amsmath,amscd,amsthm,amssymb}
\usepackage{color}

\def\Rn{{\mathbb{R}^n}}

\def\L1loc{L^{1}_{\rm loc}(\Rn)}

\newtheorem{theorem}{Theorem}{}
\newtheorem{corollary}{Corollary}{}
\newtheorem{definition}{Definition}{}
\newtheorem{remark}{Remark}{}
{}
\newtheorem{proposition}{Proposition}{}

\numberwithin{equation}{section}

\begin{document}

\begin{center}
\Large \bf Different approach to the decomposition theory of $HM^p_{q,{\Delta_{\nu}}}$ Hardy-Morrey
spaces
\end{center}

\

\centerline{\large Cansu Keskin$^{a,}$
	\footnote{ Corresponding author.
	 }, }

\

\centerline{$^{a}$ Dumlupinar University, Department of Mathematics, Kutahya, Turkey}
\centerline{e-mail: cansu.keskin@dpu.edu.tr}
\

\begin{abstract}
The Hardy-Morrey spaces related to Laplace-Bessel differential equations are introduced in terms of
maximal functions. The atomic decomposition theory
which has the same cancellation
properties of the $H^{p}_{\Delta_{\nu}}(\mathbb{R}^{n}_{+})$ Hardy spaces
is established.
\end{abstract}

\section{Introduction}

Lebesgue spaces $L^p$ and Hardy spaces $H^p$ play an important
role in function spaces theory and in harmonic analysis as well.
So many people have shown that Hardy spaces  $H^p(\mathbb{R}^n)$ $(0<p\leq \infty)$ can
be more informative than Lebesgue spaces $L^p(\mathbb{R}^n)$ when we investigate the boundedness of
some singular integral operators. For example, the classical Riesz transforms are bounded from $H^1(\mathbb{R}^n)$
to $L^1(\mathbb{R}^n)$, although they are not bounded on $L^1(\mathbb{R}^n)$.
We can also give an equivalent definition for Hardy spaces
by using the decomposition theory. This definition states that any elements of Hardy
spaces can be represented as the series of atoms.
One of the advantages of the decompositions in Hardy spaces is that we can prove
the boundedness of some operators can be verified only for the collection of atoms.

It is worthwhile to mention that the development of a theory of Hardy spaces (or other generalizations as
Hardy-Morrey) related to some classes of differential operators has attracted a great
interest. It can be developed to other function spaces. Some results in this direction can be encountered in
\cite{2}, the decomposition of Hardy-Morrey spaces,
the decomposition of Hardy spaces with variable exponent in \cite{3}, and the atomic
decomposition of Morrey spaces \cite{1}. Motivated by these advantages that Hardy spaces
enjoy, in the present paper, we investigate the atomic decomposition for
$HM^p_{q,{\Delta_{\nu}}}$ Hardy-Morrey spaces.

Hardy-Morrey spaces are important due to the connections with $L^p$, Hardy space, $BMO^{-1}$ and etc.
The maximal characterization, the atomic decomposition of this space, which have the same cancellation properties
as the classical Hardy space, were studied by Jia and Wang \cite{2,7}. These results extend those of E. M. Stein for the Hardy space $H^p(\mathbb{R}^n)(p\leq 1)$, of A. L. Mazzucato for the Besov–Morrey space \cite{10} $\mathcal{N}_{p,q,r}^s$
and of M. E. Taylor \cite{9} for the Morrey space $\mathcal{M}^p_q(q\geq 1)$.

The atomic decomposition is one of the remarkable
features for Hardy type spaces. There are many attempts of obtaining these decompositions \cite{2,7,4}. A study of the classical Hardy-Morrey spaces by using the maximal function approach is given in \cite{2} and some applications of these spaces are given in \cite{7}. Also weighted version is investigated by Ho in \cite{4}.
An interesting problem is whether we can give an atomic characterization
of Hardy-Morrey spaces associated with
Laplace-Bessel
operator $\Delta_{\nu}$ defined as follows
\begin{equation*}\label{e1.F1}
\Delta_{\nu}= \sum\limits_{i=1}^{n} \frac{\partial ^2}{\partial
x_i^2} + \frac{\nu}{x_n}\frac{\partial}{\partial x_n}, ~~~ \nu> 0.
\end{equation*}

Before finishing
we would like to point out that in this paper, maximal function characterization and atomic decomposition
are studied in the Laplace-Bessel setting of Hardy-Morrey spaces.
We introduce some new spaces called $HM^p_{q,{\Delta_{\nu}}}$ Hardy-Morrey spaces
which generalize the Morrey spaces $M^p_{q}(q>1)$ and Hardy spaces $H^p_{\Delta_{\nu}} (p<1)$ \cite{6}.

The structure of the remaining part of the present paper is as follows:
First, we recall some informations for both Hardy and Morrey spaces in Section 2. In Section 3, the maximal function characterization
is given. Finally, as applications of the atomic characterization we will obtain some estimates of the singular integral operators associated with the generalized shift operator.
We apply all our results to the Laplace-Bessel
setting, obtaining boundedness properties for $0<q\leq p<\infty$ of
the singular integral operators appearing. We also derive the characterization of the Laplace Bessel-Hardy
spaces in terms of the Riesz-Bessel transforms, extending the results in \cite{5} to a wider class of
$HM^p_{q,{\Delta_{\nu}}}$ Hardy-Morrey spaces.

\section{Hardy-Morrey spaces for $\triangle_{\nu}$ } 

The classical Hardy spaces $H^p$ are defined by \cite{6},
\begin{equation*}
H^p=\bigg\{f\in\mathcal{S}^{'}\setminus \mathcal{P}:||f||_{H^p}=\bigg|\bigg|\sup\limits_{t>0}|\phi_t\ast f|\bigg|\bigg|_{L^p}<\infty\bigg\}.
\end{equation*}
Here, $\mathcal{S}^{'}$ is the topological dual of $\mathcal{S}$ Schwartz space is the collection of all tempered distributions on
$\mathbb{R}^{n}$ equipped with the strong topology, and $\phi_t(x)=t^{-n}\phi(t^{-1}x)$ for $\phi\in \mathcal{S}$ with
$\int\phi(x)dx=1$.

The classical Morrey space was first introduced by Morrey in \cite{8}, since then a large number of investigations have been given to them by mathematicians. Recently, many authors established the boundedness
of some type operators associated with Morrey spaces.

\begin{definition}
For $p$ and $q$ satisfying $0<q\leq p<\infty$, the homogeneous Morrey spaces $M^p_q$, are defined as
\begin{equation*}
M^p_q=\{f\in L^q_{\text{loc}}: ||f||_{M^p_q}=\sup\limits_{x\in\mathbb{R}^n, R>0}|B(x,R)|^{1/p-1/q}
||f||_{L^q(B(x,R))}<\infty\}
\end{equation*}
where $B(x,R)$ is the closed ball of $\mathbb{R}^n$ with center $x$ and radius $R$.
\end{definition}
Let $j\in\mathbb{Z}$, $k\in\mathbb{Z}^n$. The set
$$Q_{jk}=\{x\in\mathbb{R}^n:2^{-j}k_i\leq x_i\leq 2^{-j}(k_i+1), i=1,\ldots,n\}$$
is a called a dyadic cube.

\begin{remark}
Note that
\begin{equation*}
||f||_{M^p_q}\approx \sup\limits_{J:\text{dyadic}}|J|^{1/p-1/q}||f||_{L^q(J)}.
\end{equation*}
\end{remark}

Based on the above notions and notations, we can state
the Hardy-Morrey spaces in Laplace-Bessel setting.

\begin{definition}
We say that a function $f\in\mathcal{S}_{+}\setminus \mathcal{P}$ belongs to the Hardy-Morrey space $HM^p_{q,{\Delta_{\nu}}}$ for $0<q\leq p<\infty$,
if
\begin{equation*}
||f||_{HM^p_{q,{\Delta_{\nu}}}}=||\sup\limits_{t>0}|\phi(t)\otimes f|||_{M^p_{q,\nu}}<\infty.
\end{equation*}
Here $\phi\in\mathcal{S}_{+}$($\mathcal{S}_{+}(\mathbb{R}^n_+)$ be the space of functions which are the restrictions
to $\mathbb{R}^{n}_{+}$ of the test functions of the Schwartz that are even with respect to $x_{n}$,
decreasing sufficiently rapidly at infinity) satisfies $\int\phi(x)x_{n}^{\nu}dx=1$ and $\mathcal{P}$ is the set of polynomials.
\end{definition}

The $HM^p_{q,{\Delta_{\nu}}}$ Hardy-Morrey spaces cover Hardy spaces. From the definition,
$H^p_{\Delta_{\nu}}=HM^p_{p,{\Delta_{\nu}}}\subset HM^p_{q,{\Delta_{\nu}}}$ for $q\leq p<\infty$.
Thus, the Hardy space $H^{p}_{\Delta_{\nu}}(\mathbb{R}^n_+)$ is the space of those
$f\in \mathcal{S}'_{+}(\mathbb{R}^n_+)$ (a function called bounded tempered distribution on $\mathbb{R}^{n}_{+})$
for which maximal function belongs to $L^p_{\nu}(\mathbb{R}^n_+)$, and defined by \cite{5}
\begin{equation*}
H^p_{\Delta_{\nu}}=\{f\in\mathcal{S}'_+\setminus \mathcal{P}:||f||_{H^p_{\Delta_{\nu}}}=||\sup\limits_{t>0}|\phi(t)\otimes f|||_{L^p_{\nu}}<\infty\}.
\end{equation*}
Here, the generalized convolution($\Delta_{\nu}$-convolution) is defined by
$$
(f \otimes g)(x)=\int_{\mathbb{R}^{n}_{+}}f(y) \; T^y g(x)\, y_n^{\nu}dy,
$$
and $T^y$ denotes the generalized shift operator given by \cite{5}.

\section{Maximal Characterization of $HM^p_{q,{\Delta_{\nu}}} $ spaces}

The Hardy-Morrey spaces are defined in this section by
using the maximal function characterizations.
One of the remarkable and fundamental properties of Hardy type spaces are the equivalence of definitions of Hardy type
spaces by different maximal functions such as the nontangential maximal function characterizations and the grand maximal function characterizations. The equivalence of these maximal function characterizations are established in this section.

We review some needed background on Hardy-Morrey spaces.
For any $N \in \mathbb{Z}^{+}$, define
$$
p_{N,\nu}(\varphi) = \sup\limits_{x\in \mathbb{R}^{n}_{+}} \big(1+|x|\big)^N \sum\limits_{|\alpha|\le N} \big| D_{x'}^{\alpha'} B_{n}^{\alpha_n} \varphi(x) \big| ~ x_n^{\nu} dx
$$
and
\begin{equation*} \label{kutahya1}
\mathfrak{F}_{N,\nu} = \Big\{ \varphi \in \mathcal{S}_{+}(\mathbb{R}^n_+) : p_{N,\nu}(\varphi) \le 1  \Big\}.
\end{equation*}
For any  positive $t$ and $\varphi \in \mathcal{S}_{+}(\mathbb{R}^{n}_{+})$, we write $\varphi_t(x)
=t^{-n-\nu}\varphi(x/t)$.

Let $\varphi \in \mathcal{S}_{+}(\mathbb{R}^n_+)$ with $\int_{\mathbb{R}^{n}_{+}} \varphi(x)~x_{n}^{\nu} dx=1$ and any $f \in \mathcal{S}'_{+}(\mathbb{R}^n_+)$, the maximal function associated with the Laplace-Bessel differential operator is defined by  V.S. Guliyev in \cite{11,12,13}
\begin{align*}
\mathcal{M}_{\nu}f(x)=\sup\limits_{r>0}\frac{1}{|E(0,r)|_{\nu}} \int_{E(0,r)}T^y|f(x)| y_n^{\nu} dy.
\end{align*}
Note that $\mathcal{M}_{\nu}$ is the Hardy-Littlewood maximal operator associated
with the Laplace-Bessel differential operator.

The grand maximal function is defined by setting
\begin{align*}
\mathcal{M}_{N,\nu}f(x)=\sup\{|(f\otimes\varphi_{t} )(x)|:t>0,~~~\varphi\in \Im_{N,\nu}\}
\end{align*}
where we choose and fix a large integer $N$.

We define the "nontangential" version of $\mathcal{M}_{\nu}$, given by
\begin{equation*}\label{kut01}
\mathcal{M}_{\varphi,\nu}^{*}f(x) = \sup\limits_{y \in \Gamma_t(x)} \left| (f \otimes \varphi_t)(y)\right|,
\end{equation*}
and also the "radial(vertical)" maximal functions defined by
\begin{align*}
\mathcal{M}^{0}_{\varphi,\nu}f(x) = \sup\limits_{0<t<\infty}\left| (f \otimes \varphi_t)(x)\right|,
\end{align*}
Here $\Gamma_t(x) = \{y\in \mathbb{R}^{n}_{+}: |x-y| <t \}$ and $\varphi_t(y)=t^{-n-\nu}\varphi\big(\frac{y}{t}\big)$.

Beside this, the $\varphi$-tangential maximal function
\begin{align*}
\mathcal{M}_{\varphi,\nu}^{\lambda}f(x)=\sup\limits_{y\in \mathbb{R}^n_{+},t>0}
|(f\otimes{\varphi_t})(y)|\bigg(\frac{t}{|x-y|+t}\bigg)^{\lambda}.
\end{align*}

Let
\begin{equation*}
\begin{array}{ll}
P^{\nu}(x,t)=P_{t}^{\nu}(x)=\dfrac{C_{k,\nu}t}{(t^2+|x|^2)^{\frac{n+k+\nu}{2}}}
,\quad C_{k,\nu}=\dfrac{2^{\frac{n+\nu}{2}}\Gamma\big(\frac{n+k+\nu}{2}\big)}{\Gamma\big({\frac{k}{2}}\big)}.
\end{array}
\end{equation*}
be the $\Delta_{\nu}$-Poisson type kernel.

\begin{definition}
The Hardy-Morrey space $HM^p_{q,{\Delta_{\nu}}}(\mathbb{R}^n_+)$ is the space of those
$f\in \mathcal{S}'_{+}(\mathbb{R}^n_+)$ for which $\mathcal{M}_{P^{\nu}}f$,
the $\Delta_{\nu}$-Poisson type maximal function and
$\mathcal{M}_{N,\nu}f$,
the grand maximal function of $f$ all
belongs to $M^p_{q,\nu}$, and define
$$
\| f \|_{HM^p_{q,{\Delta_{\nu}}}}= \| \mathcal{M}_{P^{\nu}}f \|_{M^p_{q,\nu}}=\| \mathcal{M}_{N,\nu}f \|_{M^p_{q,\nu}},
$$
where
$$
\mathcal{M}_{P^{\nu}}f(x)= \sup\limits_{t>0} \Big| \big( P^{\nu}_{t} \otimes f \big)(x)\Big|.
$$
\end{definition}

We prove one of the fundamental theorem for Hardy–Morrey spaces. It is an extension for the corresponding theorem for Hardy spaces [ See \cite{5}, Theorem 5].

\begin{theorem}\label{c}
Let $0<q\leq 1$ and $q\leq p<\infty$. Then the following statements are valid:
\begin{enumerate}
\item[(i)] There exists a Schwartz function $\varphi$ with
$\int\limits_{\mathbb{R}^n_{+}}\varphi(x)x_n^{\nu}dx=1$ and a constant $C_1$ such that
\begin{align}\label{t1}
\|\mathcal{M}_{\varphi,\nu}^{0}f\|_{M^p_{q,\nu}}\leq C_1\|f\|_{HM^p_{q,{\Delta_{\nu}}}}
\end{align}
for all $f\in \mathcal{S}'_{+}(\mathbb{R}^n_{+})$.\\
\item[(ii)]  For every $\varphi\in \mathcal{S}_{+}(\mathbb{R}^n_{+})$
there exists a constant $C_2$ such that
\begin{align}\label{4.5}
\|\mathcal{M}_{\varphi,\nu}^{*}f\|_{M^p_{q,\nu}}\leq C_2\|\mathcal{M}_{\varphi,\nu}^{0}f\|_{M^p_{q,\nu}}
\end{align}
for all $f\in \mathcal{S}'_{+}(\mathbb{R}^n_{+})$. \\
\item[(iii)]  For every $\varphi\in \mathcal{S}_{+}(\mathbb{R}^n_{+})$ and $\lambda>\frac{n+\nu}{p}$
there exists a constant $C_3$ such that
\begin{eqnarray}\label{4.6}
\|\mathcal{M}_{\varphi,\nu}^{\lambda}f\|_{M^p_{q,\nu}}\leq C_{3}\|\mathcal{M}_{\varphi,\nu}^{*}f\|_{M^p_{q,\nu}}
\end{eqnarray}
for all $f\in \mathcal{S}'_{+}(\mathbb{R}^n_{+})$. \\
\item[(iv)]  For every $\lambda\in \mathbb{Z}^{+}$ and $\varphi\in \mathcal{S}_{+}(\mathbb{R}^n_{+})$ with
$\int\limits_{\mathbb{R}^n_{+}}\varphi(x)x_n^{\nu}dx=1$
there exists a constant $C_4$ such that if $N\geq\lambda+1$ we have
\begin{eqnarray}\label{4.7}
\|\mathcal{M}_{{N,\nu}}f\|_{M^p_{q,\nu}}\leq C_4\|\mathcal{M}_{\varphi,\nu}^{\lambda}f\|_{M^p_{q,\nu}}
\end{eqnarray}
for all $f\in \mathcal{S}'_{+}(\mathbb{R}^n_{+})$. \\
\item[(v)]  For every positive integer $N$ there exists a constant $C_5$
such that every tempered distribution $f$ with $\|\mathcal{M}_{{N,\nu}}f\|_{M^p_{q,\nu}}<\infty$
is a bounded distribution and satisfies
\end{enumerate}
\begin{eqnarray}\label{t9}
\|f\|_{HM^p_{q,{\Delta_{\nu}}}}\leq C_5\|\mathcal{M}_{{N,\nu}}f\|_{M^p_{q,\nu}}
\end{eqnarray}
that is, lies in the Hardy-Morrey space $HM^p_{q,{\Delta_{\nu}}}(\mathbb{R}^n_{+})$.
\end{theorem}

We deduce that for $f\in HM^p_{q,{\Delta_{\nu}}}(\mathbb{R}^n_{+})$, the inequality
(\ref{t9}) can be reversed, and therefore for any Schwartz function $\varphi$
with  $\int\limits_{\mathbb{R}^n_{+}}\varphi(x)x_n^{\nu}dx=1$, we have
\begin{eqnarray*}
\|\mathcal{M}_{\varphi,\nu}^{*}f\|_{M^p_{q,\nu}}\leq C_6 \|f\|_{HM^p_{q,{\Delta_{\nu}}}}.
\end{eqnarray*}

We remark the pointwise following inequalities for all $x\in \mathbb{R}^n_+$
\begin{eqnarray*}
\mathcal{M}_{{N,\nu}}f(x)\leq \mathcal{M}_{\nu}f(x)\leq
\mathcal{M}_{\varphi,\nu}^{0}f\leq
C\mathcal{M}_{\varphi,\nu}^{*}f(x)\leq C\mathcal{M}_{{\varphi,\nu}}^{\lambda}f(x).
\end{eqnarray*}
It is a fact that all the maximal functions have comparable
quasinorms for all $0<q\leq 1$ and $q\leq p<\infty$.
This presents a variety of characterizations for Hardy-Morrey spaces.

\begin{proof}
The proof of the theorem is similar to [6, Theorem 5].
For a complete and detail account of the proof, the reader
is referred to [6, Theorem 5].

{\it Part (i)}. By the proof of [6, Theorem 5(i)], we have the pointwise inequalities
\begin{align}\label{4.9}
\mathcal{M}_{\varphi,\nu}^{0}f(x)
\leq C_1 \mathcal{M}_{P^{\nu}}f(x),
\end{align}
$\forall x\in\mathbb{R}^n_+$. Thus, (\ref{t1}) follows from
(\ref{4.9}).

{\it Part (ii)}. Similar to the proof [6, Theorem 5(ii)], we present the proof for $\mathcal{M}_{\varphi,\nu}^{*}f$ only as the general case follows similarly. To this end, let us now set
\begin{align*}
\mathcal{M}_{\varphi,\nu,\varepsilon,N}^{*}f=\sup_{\{y \in \mathbb{R}^n_{+}:|x-y|<t<1/\varepsilon\}}
\left|\left(f\otimes\varphi_t\right)(y)\right|\left(\dfrac{t}{t+\varepsilon}\right)^N
\left(1+\varepsilon\left|y\right|\right)^{-N}.
\end{align*}
where $0<\varepsilon<1$. Firstly, we deduce that if $\mathcal{M}_{\varphi,\nu}^{0}f\in M^p_{q,\nu}$,
then $\mathcal{M}_{\varphi,\nu}^{*}f\in
M^p_{q,\nu}\cap L^\infty_{\nu}(\mathbb{R}^n_{+})$.
From [6, p. 14], we obtain
\begin{align}\label{4.10}
\mathcal{M}_{\varphi,\nu,\varepsilon,N}^{*}
&\leq \frac{C}{\left(1+\varepsilon\left|x\right|\right)^{N-\ell}}\nonumber\\&
\leq C_{\epsilon}\bigg(\chi_{E(0,1)}(x)+\sum\limits_{k=1}^{\infty}
2^{-k(N-\ell)}\chi_{E(0,2^k)\setminus E(0,2^{k-1})}(x)\bigg).
\end{align}
Since the characteristic function is in $ M^p_{q,\nu}$, we
may write
\begin{align}\label{4.11}
||\chi_{E(0,2^k)\setminus E(0,2^{k-1})}||_{M^p_{q,\nu}}
\leq C\max(||\chi_{E(0,2^{k+1})}||_{L^{p}_{\nu}(\mathbb{R}^n_{+})}
, ||\chi_{E(0,2^k)}||_{L^{p}_{\nu}(\mathbb{R}^n_{+})}).
\end{align}
 By the inequalities (\ref{4.10}) and (\ref{4.11}), we further conclude that, for any $N$ is large enough, $\mathcal{M}_{\varphi,\nu}^{*}(f)^{\varepsilon,N} \in
M^p_{q,\nu}\cap L^\infty_{\nu}(\mathbb{R}^n_{+})$.

Next, we recall two functions used in [6, p. 14-15]. For any $M>\frac{n+\nu}{p}$, we introduce functions
\begin{align*}
V_{\varepsilon,N,M,\nu}(x)=\sup_{y\in \mathbb{R}^n_{+},t<1/\varepsilon}\left|\left(f\otimes\varphi_t\right)(y)\right
|\left(\dfrac{t}{t+\varepsilon}\right)^N\left(1+\varepsilon\left|y\right|\right)
^{-N}\left(\dfrac{t}{\left|x-y\right|+t}\right)^M
\end{align*}
and
\begin{align*}
U_{\varepsilon,N,\nu}(x)=\sup_{\{y \in \mathbb{R}^n_{+}:|x-y|<t<1/\varepsilon\}}t\left|\nabla_{\nu}
\left(f\otimes\varphi_t\right)(y)\right|\left(\dfrac{t}{t+\varepsilon}\right)^N\left(1+\varepsilon\left|y\right|\right)^{-N}
\end{align*}
where $\nabla_{\nu}=\big( \partial_1, \ldots, \partial_{n-1}, B_n \big)$.
From [6, p. 16], it is easy to see that
\begin{align*}
\left|\left(f\otimes\varphi_t\right)(y)\right
|\left(\dfrac{t}{t+\varepsilon}\right)^N
\left(1+\varepsilon\left|y\right|\right)
^{-N}\leq
\left(\dfrac{t}{\left|x-y\right|+t}\right)^M
\{\mathcal{M}_{\nu}[\mathcal{M}_{\varphi,\nu,\varepsilon,N}^{*}f]
^{s}\}^{1/s}(x)
\end{align*}
where $0<s<\infty$. By the properties of the $\triangle_{\nu}$-maximal function (see \cite{11,12}), we also have
\begin{align}\label{4.13}
||V_{\varepsilon,N,M,\nu}||_{M^p_{q,\nu}}&\leq C
||\mathcal{M}_{\nu}[\mathcal{M}_{\varphi,\nu,\varepsilon,N}^{*}f]^s||
^{1/s}_{M^p_{q/s,\nu}}\nonumber\\&
\leq C||[\mathcal{M}_{\varphi,\nu,\varepsilon,N}^{*}f]^s||
^{1/s}_{M^p_{q/s,\nu}}
\leq C||\mathcal{M}_{\varphi,\nu,\varepsilon,N}^{*}f||_{M^p_{q,\nu}}
\end{align}
for some $C>0$ independent of $f$. Moreover, according to
[6, (15)], we obtain the pointwise inequality
\begin{align}\label{4.14}
U_{\varepsilon,N,\nu}(x)\leq CV_{\varepsilon,N,M,\nu}(x).
\end{align}
Hence (\ref{4.13}) and (\ref{4.14}) yield
\begin{align}\label{4.15}
||U_{\varepsilon,N,\nu}||_{M^p_{q,\nu}}\leq C_0 ||\mathcal{M}_{\varphi,\nu,\varepsilon,N}^{*}f||_{M^p_{q,\nu}}
\end{align}
where $C_0$ is a constant to be determined in [6, p. 16].
Now, if we set $E_{\varepsilon}=\{x:U_{\varepsilon,N,\nu}(x)\leq 2C_{0}\mathcal{M}_{\varphi,\nu,\varepsilon,N}^{*}f(x)\}$,
where $C_0$ is a constant to be determined in (\ref{4.15}), we may write
\begin{eqnarray*}
||\chi_{(E_{\varepsilon})^{c}}
\mathcal{M}_{\varphi,\nu,\varepsilon,N}^{*}f||_{M^p_{q,\nu}}
\leq \frac{1}{2C_{0}}||\chi_{(E_{\varepsilon})^{c}}
U_{\varepsilon,N,\nu}||_{M^p_{q,\nu}}
\leq\frac{1}{2C_0}||U_{\varepsilon,N,\nu}||_{M^p_{q,\nu}}.
\end{eqnarray*}
It follows from (\ref{4.15}),
\begin{eqnarray}\label{4.16}
||\chi_{(E_{\varepsilon})^{c}}
\mathcal{M}_{\varphi,\nu,\varepsilon,N}^{*}f||_{M^p_{q,\nu}}
\leq \frac{1}{2}
||\mathcal{M}_{\varphi,\nu,\varepsilon,N}^{*}f||_{M^p_{q,\nu}}.
\end{eqnarray}
Moreover, by  [6, p. 16], we obtain
\begin{eqnarray*}
\mathcal{M}_{\varphi,\nu,\varepsilon,N}^{*}f\leq C\{\mathcal{M}_{\nu}[\mathcal{M}_{\varphi,\nu}^{0}f]^{s}(x)\}^{1/s}, ~\text{if}~
x\in E_{\varepsilon}
\end{eqnarray*}
where $0<s<\infty$.

Similar to the proof of  (\ref{4.13}), we have
\begin{eqnarray}\label{4.17}
||\chi_{(E_{\varepsilon})}
\mathcal{M}_{\varphi,\nu,\varepsilon,N}^{*}f||_{M^p_{q,\nu}}
\leq C
||\mathcal{M}_{\varphi,\nu}^0 f||_{M^p_{q,\nu}}.
\end{eqnarray}
The inequalities (\ref{4.16}) and (\ref{4.17}) yield the following estimate
\begin{eqnarray*}
||
\mathcal{M}_{\varphi,\nu,\varepsilon,N}^{*}f||_{M^p_{q,\nu}}
\leq C
||\mathcal{M}_{\varphi,\nu}^0 f||_{M^p_{q,\nu}}+\frac{1}{2}
||\mathcal{M}_{\varphi,\nu,\varepsilon,N}^{*}f||_{M^p_{q,\nu}}.
\end{eqnarray*}
Since $\mathcal{M}_{\varphi,\nu,\varepsilon,N}^{*}f\in M^p_{q,\nu}$, by applying Lebesgue monotone theorem, we obtain
\begin{eqnarray*}
||
\mathcal{M}_{\varphi,\nu}^{*}f||_{M^p_{q,\nu}}
\leq C||\mathcal{M}_{\varphi,\nu}^0 f||_{M^p_{q,\nu}}
\end{eqnarray*}
for some $C>0$ depending on $N$. This assure that
\begin{eqnarray*}
||\mathcal{M}_{\varphi,\nu}^0 f||_{M^p_{q,\nu}}<\infty
\Rightarrow ||
\mathcal{M}_{\varphi,\nu}^{*}f||_{M^p_{q,\nu}}<\infty.
\end{eqnarray*}
With this assertion, (\ref{4.5}) can be established by repeating the above arguments with $U_{\varepsilon,N,\nu}$ and $V_{\varepsilon,N,M,\nu}$ replaced by two new auxiliary functions
\begin{align*}
U_{\varepsilon,N,\nu}(x)=\sup_{\{y \in \mathbb{R}^n_{+}:|x-y|<t<1/\varepsilon\}}t\left|\nabla_{\nu}
\left(f\otimes\varphi_t\right)(y)\right|\left(\dfrac{t}{t+\varepsilon}\right)^N\left(1+\varepsilon\left|y\right|\right)^{-N}
\end{align*}
and
\begin{align*}
V_{\varepsilon,N,M,\nu}(x)=\sup_{y\in \mathbb{R}^n_{+},t<1/\varepsilon}\left|\left(f\otimes\varphi_t\right)(y)\right
|\left(\dfrac{t}{t+\varepsilon}\right)^N\left(1+\varepsilon\left|y\right|\right)
^{-N}\left(\dfrac{t}{\left|x-y\right|+t}\right)^M.
\end{align*}
For a detail account the reader is referred to [6, pp. 14-17].

{\it Part (iii)}. Again return to the proof [6, Theorem 6(iii)],
we have the following inequality
\begin{eqnarray*}
\{\mathcal{M}_{\varphi,\nu}^{\lambda}f(x)\}^{q}\leq C_3 \mathcal{M}_{\nu}[\mathcal{M}_{\varphi,\nu}^{*}f]^{q}(x).
\end{eqnarray*}
Then, by the above inequality with $q<p$ and the boundedness properties of $\triangle_{\nu}$-maximal functions in
$M^p_{q,\nu}$ Morrey spaces give (\ref{4.6}) whenever $\lambda>\frac{n+\nu}{p}$.

{\it Part (iv)}. The proof of [6, Theorem 6(iv)] asserts that
\begin{eqnarray*}
\mathcal{M}_{{N,\nu}}f\leq C_{4}\mathcal{M}_{\varphi,\nu}^{\lambda}f.
\end{eqnarray*}
Clearly, (\ref{4.7}) follows.

{\it Part (v)}. Let $f\in \mathcal{S}'(\mathbb{R}^{n}_{+})$ satisfy
$||\mathcal{M}_{{N,\nu}}f||_{M^p_{q,\nu}}<\infty$ for some $N\in\mathbb{N}$. For any $\theta\in \mathcal{S}(\mathbb{R}^{n}_{+})$, we obtain a  fixed constant $c_0>0$ such that $c_0 \theta\in \Im_{N,\nu}$. Hence,
$\mathcal{M}_{c_0\theta,\nu}^{*}f\leq \mathcal{M}_{{N,\nu}}f$.

By the proof of [6, Theorem 6(v)], we obtain the following stronger estimate
\begin{eqnarray*}
\mathcal{M}_{P^{\nu}}f(x)=\sup\limits_{t>0}|(P^{\nu}_t\otimes f)(x)|\leq C_5\mathcal{M}_{{N,\nu}}f(x),
\end{eqnarray*}
hence, (\ref{t9}) follows.

Finally, it follows that the quasinorm $\|\mathcal{M}_{P^{\nu}}f(x)\|_{M^p_{q,\nu}}$
is also equivalent to $\|f\|_{HM^p_{q,{\Delta_{\nu}}}}$.
\end{proof}

\begin{corollary}
The Hardy-Morrey space $HM^p_{q,{\Delta_{\nu}}}$ is complete in the metric $d(f,g)=||f-g||_{HM^p_{q,{\Delta_{\nu}}}}^q$ for $q\leq 1$.
\end{corollary}

We note that if $f\in HM^p_{q,{\Delta_{\nu}}}$ and $\varphi\in
\mathcal{S}(\mathbb{R}^n_+)$, then
\begin{eqnarray}\label{2.11}
|f\otimes \varphi(x)|^q\leq \frac{C}{|B(x,1)|_{\nu}}
\int_{B(x,1)}\mathcal{M}_{{N,\nu}}^{*}f(y)y_n^{\nu}dy.
\end{eqnarray}
Hence if $f_n\to f$ in $f\in HM^p_{q,{\Delta_{\nu}}}$, then $f_n\to f$ in $\varphi\in\mathcal{S}'$. So, the embedding is valid.

\section{Atomic decomposition}
The main result of this paper is presented
in this section. The atomic decomposition for the Hardy-Morrey spaces is established.
We next define a $(p,q,s)$-atom associated to the operator
$\triangle_{\nu}$.

\begin{definition}
Let $0<q\leq 1$, $q\leq p<\infty$ and $s\in\mathbb{N}\cup\{0\}$. A $(p,q,s)$-atom $a(x)$
is a function which satisfies the following properties:
\begin{enumerate}
\item[(i)] $a$ be supported on a cube $3Q$, namely, $\text{supp}~~a\subset 3Q$,
\item[(ii)] $\|a(x)\|_{L^{\infty}_{\nu}}\leq |Q|_{\nu}^{-\frac{1}{p}}$,
\item[(iii)] $\int_{\mathbb{R}^n_+}a(x)x^{\alpha}x_n^{\nu}dx=0$  for all $s\geq [(n+k+\nu)\big(\frac{1}{q}-1\big)]$ with
$|\alpha|\leq s,$.
\end{enumerate}
Here $3Q$ is the cube concentric with $Q$
of side-length $3\ell(Q)$.
\end{definition}

The so called atomic decomposition theorem in $HM^p_{q,{\Delta_{\nu}}}$ space is as follows:
\begin{theorem}\label{t2.1}
Let $0<q\leq 1$, $q\leq p<\infty$. Then for any $f\in HM^p_{q,{\Delta_{\nu}}}$, there exist a sequence of
scalars $\lambda_j$, a sequence of $(p,q,s)$-atoms $a_j=\{a_j\}$ for
$HM^p_{q,{\Delta_{\nu}}}$ such that
\begin{align}\label{1.3}
f=\sum\limits_{j=1}^{\infty}\lambda_{j}a_{j}
\end{align}
with
\begin{align*}
||\lambda||_{p,q,\nu}=
\bigg\{\sup\limits_{Q}\bigg(\frac{1}{|Q|_{\nu}}\bigg)^{1-q/p}
\sum\limits_{j=1}^{\infty}|Q|^{1-q/p}_{\nu}|\lambda_{j}|^{q}\bigg\}
^{1/q}<\infty.
\end{align*}
The sum converges in $\mathcal{S}^{'}_{+}\setminus \mathcal{P}$
and $f\in HM^p_{q,{\Delta_{\nu}}}$ with $||f||_{HM^p_{q,{\Delta_{\nu}}}}\leq C||\lambda||_{p,q,\nu}$ for some $C=C(n,p,q,\nu)$.
Conversely, every function $f\in HM^p_{q,{\Delta_{\nu}}}$  has the atomic decomposition (\ref{1.3})
in $\mathcal{S}^{'}_{+}\setminus \mathcal{P}$, here $a_j$'s are $(p,q,s)$-atoms and
$\lambda$ satisfies $||\lambda||_{p,q,\nu}\leq C||f||_{HM^p_{q,{\Delta_{\nu}}}}$ for some $C=C(n,p,q,\nu)$.
\end{theorem}

Before proving the decomposition theorem, we state a well known preliminary supporting result for the atomic decomposition.

\begin{proposition}\label{p5.4}
Suppose that $f\in HM^p_{q,{\Delta_{\nu}}}$ with
$0<q\leq 1$, $q\leq p<\infty$ and $\delta>0$. Then for any
$f \in \mathcal{S}_{+}(\mathbb{R}^n_+)$, there exist
$g \in \mathcal{S}_{+}'(\mathbb{R}^n_+)$, a collection of cubes
$\{Q\}_{k\in\mathbb{N}}$ and a family of smooth functions
with compact supports $\{\xi_k\}$ such that

\begin{itemize}
  \item [(i)] $f=g+b$ where $b=\sum\limits_{k\in\mathbb{N}}b_k$,
  \item [(ii)] The cubes $\{Q_k\}_{k\in\mathbb{N}}$
  which are mutually disjoint, satisfy
 \begin{align*}
 \bigcup\limits_{k\in\mathbb{N}}Q_k=
 \{x\in\mathbb{R}^n_{+}:(\mathcal{M}_{N,\nu}f)(x)>\delta\},
 \end{align*}
  \item [(iii)] $\text{supp}~\xi_k\subset Q_k$,
  $0\leq \xi_k\leq 1$ and $\sum\limits_{k\in\mathbb{N}}\xi_k=\chi_{ \{x\in\mathbb{R}^n_{+}:(\mathcal{M}_{N,\nu}f)(x)>\delta\}}$
  \item[(iv)] The function $g$ satisfies
   \begin{align*}
  (\mathcal{M}_{N,\nu}g)(x)&\leq (\mathcal{M}_{N,\nu}f)(x)
  \chi_{ \{x\in\mathbb{R}^n_{+}:(\mathcal{M}_{N,\nu}f)(x)\leq \delta\}}(x)\\&
  +\delta \sum\limits_{k\in\mathbb{N}}\dfrac{\ell(Q_k)^{n+\nu+s+1}}
  {(\ell(Q_k)+|x-x_k|)^{n+\nu+s+1}},
   \end{align*}
  \item [(v)] The function $b_k$ is given by $b_k=(f-c_k)\xi_k$
  where $c_k\in \mathcal{P}_s$ satisfying
     \begin{align*}
  \int_{\mathbb{R}^n_{+}}b_k(x)m(x)x_n^{\nu}dx=0,~~~\forall m\in\mathcal{P}_s,
    \end{align*}
    and
    \begin{align*}
  (\mathcal{M}_{N,\nu}b_k)(x)\leq C(\mathcal{M}_{N,\nu}f)(x)
  \chi_{Q_k}(x)+\delta \dfrac{\ell(Q_k)^{n+\nu+s+1}}
  {|x-x_k|^{n+\nu+s+1}}\chi_{\mathbb{R}^n_{+}\setminus Q_k}(x)
     \end{align*}
     for some $C>0$.
\end{itemize}
\end{proposition}

We are now ready to prove the main theorem.

{\it Proof of Our Main Theorem}.
In the rest of this section, for simplicity, we use
the notation $\mathcal{M}f=\mathcal{M}_{N,\nu}f$.

It suffices to show the atomic decomposition for $(p,\infty,s)$ atoms since $(p,\infty,s)$
atoms are $(p,q,s)$ atoms for any $1\leq q<\infty$.
According to the Proposition \ref{p5.4}, for any $\delta=2^j, j\in\mathbb{Z}$, we have distributions $g^j, b^j$ satisfying the properties (i)-(v) in the Proposition \ref{p5.4} and
$f=g^j+b^j$. We can write

\begin{align}\label{5.9}
O^j=\{x\in\mathbb{R}^n_{+}: (\mathcal{M}f)(x)>2^j\}
=\bigcup\limits_{k\in\mathbb{N}}Q^j_k
\end{align}
where $\bigcup\limits_{k\in\mathbb{N}}Q^j_k$ is the decomposition of $O^j$ given in the Proposition below.

Let $\{\xi_k^j\}$ be the family of smooth functions with respect to decomposition $O^j$ given by Proposition \ref{p5.4}, property (iii). In addition, as $\mathcal{M}f\in M^p_q$, we obtain $O^{j+1}\subset O^j$
and  for all $j\in\mathbb{N}$, $O^j$ are mutually disjoint. For any $\varphi\in\mathcal{S}(\mathbb{R}^n)$, we have a constant $c>0$ such that $c\varphi\in \mathfrak{F}_{N,\nu}$.
Proposition \ref{p5.4} gives

\begin{align*}
c|\varphi\ast g^j(x)|&\leq (\mathcal{M}g^j)(x)\\&
\leq (\mathcal{M}f)(x) \chi_{x\in \mathbb{R}^n_{+}:(\mathcal{M}f)(x)\leq 2^j}(x)
+2^j \sum\limits_{k\in\mathbb{N}}\dfrac{\ell(Q_k)^{n+\nu+s+1}}
{(\ell(Q_k)+|x-x_k^j|)^{n+\nu+s+1}}\leq C2^j
\end{align*}
for some $C>0$ where $x_k^j$ is the center of $Q_k^j$.
That is $g^j\to 0$ in $\mathcal{S}'(\mathbb{R}^n_{+})$ as
$j\to-\infty$.

Observe next that $b^j\to 0$ in $\mathcal{S}'(\mathbb{R}^n_{+})$ as $j\to\infty$. By
the item (v) of Proposition \ref{p5.4}, similar to the above proof works and (\ref{5.9}), for any $n/(n+\nu+s+1)<q$,
we obtain
\begin{align}\label{3.21}
\int_{Q}|(\mathcal{M}b^j)(x)|^q x_n^{\nu}dx
\leq \int_{O^j}|(\mathcal{M}f)(x)|^q (\mathcal{M}_{\nu}\chi_{Q})(x)x_n^{\nu}dx
\end{align}
for some $C>0$.

In fact, for any $\varphi\in\mathcal{S}(\mathbb{R}^n_{+})$, using the inequalities (\ref{2.11}) and (\ref{3.21}), we find that
\begin{align*}
|b^j\otimes \varphi(x)|^q&\leq C\frac{1}{|Q(x,1)|_{\nu}}
\int_{Q(x,1)}|\mathcal{M}_{\varphi,\nu}^{*}(b^j )(y)|^q
y_n^{\nu}dy\leq C\int_{Q(x,1)}|(\mathcal{M}b^j)(y)|^q y_n^{\nu}dy\\&
\leq C\int_{O^j}|(\mathcal{M}f)(y)|^q (\mathcal{M}_{\nu}\chi_{Q(x,1)})(y)
y_n^{\nu}dy\\&
C\int_{O^j}|(\mathcal{M}f)(y)|^q (1+|x-y|)^{-n-\nu}
y_n^{\nu}dy.
\end{align*}

By using the Hölder inequality for the pair
$L_{p/q,\nu}(\mathbb{R}^n_{+})$ and
$L_{(p/q)',\nu}(\mathbb{R}^n_{+})$, we obtain

\begin{align*}
\int_{\mathbb{R}^n_{+}}|(\mathcal{M}f)(y)|^q (1+|x-y|)^{-n-\nu}y_n^{\nu}dy&\leq
C\sum\limits_{k=0}^{\infty} 2^{-k(n+\nu)}
\int_{\mathbb{R}^n_{+}}|(\mathcal{M}f)(y)|^q
\chi_{B_k}(y)y_n^{\nu}dy\\&
\leq C \sum\limits_{k=0}^{\infty}\frac{1}{|B(x,2^k)|_{\nu}}
||(\mathcal{M}f)(y)^q||_{L_{\nu}^{p/q}(\mathbb{R}^n_{+})}
||\chi_{B(x,2^k)}||_{L_{\nu}^{(p/q)'}(\mathbb{R}^n_{+})}\\&
\leq C \sum\limits_{k=0}^{\infty}\dfrac{||\mathcal{M}f||^q_{M_{q,\nu}^p}}
{||\chi_{B(x,2^k)}||_{L_{\nu}^{p/q}(\mathbb{R}^n_{+})}}
\leq C ||\mathcal{M}f||^q_{M_{q,\nu}^p},
\end{align*}
where $B_k=B(x,2^k)\setminus B(x,2^{k-1})$ when $k\geq 1$
and $B_{0}=B(x,1)$.

In view of the fact that $O^j\downarrow \emptyset$, the dominated convergence theorem gives that
\begin{align*}
\lim\limits_{j\to \infty} |b^j\otimes \varphi(x)|^q
\leq C \lim\limits_{j\to \infty}\int_{O^j}
|(\mathcal{M}f)(y)|^q (1+|x-y|)^{-n-\nu}
y_n^{\nu}dy=0.
\end{align*}
Thus, $b^j\otimes \varphi\to 0$ pointwisely. That is,
$b^j\to 0$ in $\mathcal{S}'(\mathbb{R}^n_{+})$ when $j\to\infty$.

Now we return to the atomic decomposition. The telescoping
sum $f=\sum\limits_{k\in\mathbb{Z}}(g^{j+1}-g^j)$ converges in
the sense of distributions. By using the property $(v)$
of Proposition \ref{p5.4}, we obtain

\begin{align*}
g^{j+1}-g^j=b^{j+1}-b^j=\sum\limits_{k\in\mathbb{N}}
((f-c_k^{j+1})\xi_k^{j+1}-(f-c_k^j)\xi_k^{j})
\end{align*}
where $c_k^j\in\mathcal{P}_s$. Also, $\forall m\in\mathcal{P}_s$
\begin{align*}
\int_{\mathbb{R}^n_{+}}
(f(x)-c_k^{j}(x))m(x)\xi_k^{j}(x)x_n^{\nu}dx=0
\end{align*}
holds.
By the fact $\sum_{k\in\mathbb{N}}\xi_k^{j}=1$ on the support
$\xi_k^{j+1}$, we obtain
\begin{align}\label{3.22}
f=\sum_{j,k}A^{j}_k=\sum_{j}(g^{j+1}-g^j),
\end{align}
 where
\begin{align*}
A^{j}_k=(f-c_k^{j})\xi_k^{j}-\sum\limits_{l\in\mathbb{N}}
(f-c_l^{j+1})\xi_l^{j+1}\xi_k^{j}+\sum\limits_{l\in\mathbb{N}}
c_{k,l}\xi_l^{j+1},~~~~c_{k,l}\in\mathcal{P}_s
\end{align*}
and
\begin{align*}
a_k^j=\lambda^{-1}_{j,k}A^j_k~~~~~~~~\text{and}~~~~
\lambda_{j,k}=C2^j||\chi_{Q^j_k}||_{L^p_{\nu}(\mathbb{R}^n_{+})}.
\end{align*}
The function $a_k^j$'s are obviously $(p,\infty,s)$ atoms.

According to the definition of $Q_k^j$ and following from the fact that the family $\{Q^j_k\}_{k\in\mathbb{N}}$ has the finite intersection property, we compute that for any $0<q<\infty$

\begin{align*}
\sum\limits_{k\in\mathbb{N}}\bigg(\frac{|\lambda_{j,k}|}
{||\chi_{Q^j_k}||_{L^p_{\nu}(\mathbb{R}^n_{+})}}\bigg)^q
\chi_{Q^j_k}(x)\leq C2^{qj}\chi_{O^j}(x)
\end{align*}
for some $C>0$. Consequently,
\begin{align*}
\sum\limits_{j,k}\bigg(\frac{|\lambda_{j,k}|}
{||\chi_{Q^j_k}||_{L^p_{\nu}(\mathbb{R}^n_{+})}}\bigg)^s
\chi_{Q^j_k}(x)\leq  C\sum\limits_{j\in\mathbb{Z}}
2^{qj}\chi_{O^j}(x)\leq C(\mathcal{M}f)(x)^q.
\end{align*}

By applying the quasi-norm $||.||^{1/q}_{M_{p/q,\nu}}$
on both sides of the above inequality,
\begin{align*}
\bigg|\bigg|\sum\limits_{j,k}\bigg(\frac{|\lambda_{j,k}|}
{||\chi_{Q^j_k}||_{L^p_{\nu}(\mathbb{R}^n_{+})}}\bigg)^q
\chi_{Q^j_k}\bigg|\bigg|^{1/q}_{M_{p/q,\nu}}
\leq  C||f||_{HM^p_{q,{\Delta_{\nu}}}},~~~0<q<\infty
\end{align*}
for some $C>0$ independent of $f$.

It follows from that
\begin{align*}
||\lambda||_{p,q,\nu}\leq C||f||_{HM^p_{q,{\Delta_{\nu}}}}.
\end{align*}
We can replace (\ref{3.22}) by
$
f=\sum\limits_{j=1}^{\infty}\lambda_{j}a_{j}.
$
Thus, the proof of our main theorem is completed.

\subsection*{Acknowledgement} The author would like to thank the reviewers for valuable suggestions and corrections.

\bigskip

\end{document}